\documentclass[11pt]{article}
\usepackage{latexsym}
\usepackage{amsmath,amsthm,amsfonts,amssymb,graphicx}
\title{Van Kampen's embedding obstruction for discrete groups} \author{Mladen
  Bestvina, Michael Kapovich, and Bruce Kleiner \thanks{All three authors
  gratefully acknowledge the support by the
    National Science Foundation.}}  \date{\today}

\newtheorem{thm}{Theorem}%[section]
\newtheorem{lemma}[thm]{Lemma}
\newtheorem{cor}[thm]{Corollary}
\newtheorem{prop}[thm]{Proposition}
{}

\theoremstyle{remark}
\newtheorem{example}[thm]{Example}

\newtheorem{definition}[thm]{Definition}
\newtheorem{remark}[thm]{Remark}

\def\R{{\mathbb R}}
\def\Z{{\mathbb Z}}
  
\def\G{{\Gamma}}
\newcommand{\<}{\langle}
\renewcommand{\>}{\rangle}

\begin{document}
\maketitle
\begin{abstract} We give a lower bound to the dimension of a contractible
  manifold on which a given group can act properly discontinuously. In
  particular, we show that the $n$-fold product of nonabelian free groups
  cannot act properly discontinuously on $\R^{2n-1}$.
\end{abstract}
\def\obdim{\operatorname{obdim}}
\def\qedim{\operatorname{updim}}
\def\actdim{\operatorname{actdim}}
\def\gdim{\operatorname{gdim}}
\def\vcdim{\operatorname{vcdim}}
\def\diam{\operatorname{diam}}

\section{Introduction}

In \cite{vk:embedding} van Kampen developed an obstruction theory for
embeddings of finite $n$-complexes into $\R^{2n}$. We will briefly review van
Kampen's theory in Section \ref{section:recall}. It is natural (and
straightforward) to remove the dimension restrictions and talk about
a cohomological obstruction to embedding a complex into $\R^m$. Complexes where
this obstruction does not vanish will be called {\it $m$-obstructor
  complexes}. The precise definition will be given below (see Definition
\ref{obstructor}). For example, the utilities graph (the join of two 3-point
sets) is a 2-obstructor complex, and van Kampen proved that the $n$-fold join
of 3-point sets is a $(2n-2)$-obstructor complex.

We introduce the notion of the {\it obstructor dimension} $\obdim\G$ of a
discrete group $\G$ (see Definition \ref{obdim}). For example, when the group
is hyperbolic or CAT(0) and the boundary contains an $m$-obstructor complex,
then $\obdim\G\geq m+2$. In particular, $\obdim F_2^n=2n$ ($F_2$ is the free
group of rank 2) since the boundary of $F_2^n$ is the $n$-fold join of Cantor
sets and thus contains the complex considered by van Kampen.

The main theorem in this paper is the following.

\begin{thm} \label{1}
If $\obdim\G\geq m$ then $\G$ cannot act properly discontinuously
  on a contractible manifold of dimension $<m$.
\end{thm}

\noindent For example, $F_2^n$ cannot act properly discontinuously on $\R^{2n-1}$.

\noindent In \cite{bf:cusp} the methods of this paper are used to prove:

\begin{thm}\cite{bf:cusp} \label{2}
Let $G$ be a connected semisimple Lie group,
  $K\subset G$ a maximal compact subgroup, $G/K$ the associated contractible
  manifold (i.e., the symmetric space when the center of $G$ is finite) and
  $\G$ a lattice in $G$. If $\G$ acts properly discontinuously on a
  contractible manifold $W$, then $\dim W\geq \dim G/K$.\end{thm}

There is an application of our results to Geometric Topology.
The celebrated theorems of Whitney \cite{whitney:embeddings},
\cite{whitney:immersions} state that every $n$-manifold can be embedded in
$\R^{2n}$ and immersed (for $n>1$) into $\R^{2n-1}$. A less known theorem
of Stallings \cite{dr:stallings} asserts that every $n$-complex is homotopy
equivalent to a complex that embeds in $\R^{2n}$. It is therefore natural to
ask whether every $n$-complex can be immersed up to homotopy into
$\R^{2n-1}$.

\begin{cor}
Let $X=X_1\times X_2\times \cdots\times X_n$ be the $n$-fold product of connected graphs
$X_i$ with the first betti number 2. Then $X$ does not immerse up
to homotopy into $\R^{2n-1}$. 
\end{cor}

\begin{proof}
Suppose $X\simeq Y$ and $Y$ immerses in $\R^{2n-1}$. Then $Y$ has a thickening
(an immersed in $\R^{2n-1}$ regular neighborhood) which is an aspherical $(2n-1)$-manifold
with fundamental group $F_2^n$. But then the universal covering action
violates Theorem \ref{1}.
\end{proof}

It appears that the above comlplex is the first example 
of aspherical simplicial complex of dimension $n\ge 3$ which 
does not immerse up to homotopy into $\R^{2n-1}$. 
A detailed study of thickenings in the case when $K=S^m\cup_\alpha e^n$
has 3 cells was carried out by Cooke \cite{cooke:thickenings}. In particular,
he constructs such complexes where the minimal thickenings have arbitrarily
large codimension.

\section{Obstructor complexes}\label{section:recall}

In this section we briefly recall the work of van Kampen \cite{vk:embedding}.
At the time his paper was written, cohomology theory was still not fully
developed and many of the details were elaborated later in
\cite{shapiro:embedding} (see also \cite{wu:embedding}).  Van Kampen
constructed an $n$-complex that does not embed into $\R^{2n}$. Van Kampen's
complex is the $(n+1)$-fold join of the 3-point set, generalizing the
well-known non-planar ``utilities'' graph.  Flores \cite{flores:embedding}
showed that the $n$-skeleton of the $(2n+2)$-simplex works just as well, thus
generalizing the other standard example of a non-planar graph, namely the
complete graph on 5 vertices. Flores reduced the claim to the
Borsuk-Ulam theorem and no additional cohomological arguments were needed.

In what follows, we shall also need examples of
complexes that embed in an even dimensional Euclidean space, but not in one of
lower dimension. We will follow the standard practice and blur the distinction
between a simplicial complex and its geometric realization. All (co)homology
groups are taken with coefficients in ${\mathbb Z}_2$.

\begin{definition}\label{obstructor}
  Fix a non-negative integer $m$. A finite simplicial complex $K$ of dimension
  $\leq m$ is an {\it $m$-obstructor complex} if the following holds:
\begin{enumerate}
\item There is a collection 
$$
\Sigma=\{\{\sigma_i,\tau_i\}_{i=1}^k\}
$$
  of
  unordered pairs of disjoint simplices of $K$ with
  $\dim\sigma_i+\dim\tau_i=m$ that determine an $m$-cycle
  (over $\Z_2$) in $$\bigcup\{ \sigma\times\tau\subset K\times
    K|\sigma\cap\tau=\emptyset\}  /\Z_2$$
  where $\Z_2$ acts by $(x,y)\mapsto (y,x)$.
\item For some (any) general position map $f:K\to \R^m$ the (finite) number 
$$
\sum_{i=1}^k  |f(\sigma_i)\cap f(\tau_i)|$$ is odd.
\item For every $m$-simplex $\sigma\in K$ the number of vertices $v$ such that
  the unordered pair $\{\sigma,v\}$ is in $\Sigma$ is even.
\end{enumerate}
\end{definition}

It turns out (see \cite{shapiro:embedding}, \cite{wu:embedding}, 
\cite{Freedman-Krushkal-Teichner})  
that Van Kampen's obstruction (conditions 
(1) and (2) above) is the only 
obstruction to existence of embedding of complexes of dimension 
$\ge 3$ into $\R^{2n}$. For 2-dimensional complexes there are other 
obstructions as well, see \cite{Freedman-Krushkal-Teichner}, 
\cite{Krushkal}. 

\subsection{Discussion and basic properties}
Let $K$ be an $m$-obstructor complex. If $f$ and $f'$ are two general position maps $K\to \R^m$ choose a general
position homotopy $H$ between them. A standard argument of ``watching $H$''
shows that in the presence of item 1 (in Definition \ref{obstructor}) 
the two integers from item 2 for $f$ and
$f'$ differ by an even integer. In particular, item 2 implies that $K$ does
not embed in $\R^m$; indeed for every map $K\to \R^m$ there exist two disjoint
simplices of $K$ whose images intersect.

We will view $\Sigma$ as a subcomplex of $\bigcup\{ \sigma\times\tau\subset
K\times K|\sigma\cap\tau=\emptyset\} /\Z_2$. Then item 1 states that $\Sigma$
is a $m$-pseudomanifold over $\Z_2$, meaning that every $(m-1)$-cell is the
face of an even number of $m$-cells, and in particular we have the fundamental
class $[\Sigma]\in H_m(\Sigma)$. Similarly, the collection $\tilde \Sigma$ of
{\it ordered} pairs corresponding to the pairs in $\Sigma$ can be viewed as a
subcomplex of $\bigcup\{ \sigma\times\tau\subset K\times
K|\sigma\cap\tau=\emptyset\}$ and is an $m$-pseudomanifold. Further,
$(x,y)\mapsto (y,x)$ is the deck transformation of the natural double cover
$\tilde\Sigma\to\Sigma$. Let $\phi:\Sigma\to\R P^{\infty}$ be a classifying
map for this double cover. We note that item 2 is equivalent to the
requirement that
$$
\<\phi^*(w^m),[\Sigma]\>\neq 0\in\Z_2
$$
where $w^m\in H^m(\R P^\infty)$ is the nonzero class. Indeed, we can perturb
$f:K\to \R^m$ to a map $F=(f,g):K\to\R^m\times\R=\R^{m+1}$ so that
$F(\sigma_i)\cap F(\tau_i)=\emptyset$ for all $i$. Then we have a classifying
map $\phi:\Sigma\to\R P^m\subset \R P^\infty$ defined by 
$$
\phi(\{x,y\})=\text
{line through }F(x)\text{ and }F(y)$$
 where a point of $\R P^m$ is viewed as
the set of parallel lines in $\R^m$. Then $\<\phi^*(w^m),[\Sigma]\>$ can be
computed as the ``degree'' of $\phi$, which in turn is the number of points of
$\Sigma$ mapped to the ``vertical lines'' $pt\times \R$, i.e., the number from
item 2.

We could have defined the notion of an $m$-obstructor complex by requiring
only items 1 and 2. This definition would then be equivalent to the
requirement that $\Phi^*(w^m)\neq 0$, where $\Phi:\big( K\times K\setminus
\Delta\big) /\Z_2\to\R P^\infty$ is the classifying map ($\Delta\subset
K\times K$ is the diagonal), and this would be closer in spirit to van
Kampen's work. We impose item 3 to ensure that the Join Lemma and the Linking
Lemma below hold. A restatement of item 3 is that the projection map
$\pi:\tilde\Sigma\to K$ (say to the second coordinate) has the property that
the pullback of every $m$-cocycle evaluates trivially on the fundamental class
$[\tilde\Sigma]$.

Van Kampen's obstruction theory can be summarized in the following
proposition.

\begin{prop} \label{vk}
  Suppose that $K$ is an $m$-obstructor complex and that $W$ is a contractible
  $m$-manifold. Then for every map $F:K\to W$ there exist disjoint simplices
  $\sigma$ and $\tau$ in $K$ such that $F(\sigma)\cap F(\tau)\neq\emptyset$.
  In particular, $K$ does not embed into $W$.
\end{prop}

\begin{proof}
  The case of $W=\R^m$ was discussed above. For the general case, assume on
  the contrary that $F:K\to W$ violates the proposition. 
  Define $\phi:\Sigma\to W\times
  W\setminus\Delta/\Z_2$ by $\phi(\{x,y\})=\{F(x),F(y)\}$. The following lemma
  then implies that $\Sigma$ classifies into $\R P^{m-1}$, a contradiction.
\end{proof}

\begin{lemma}
  Suppose $W$ is a contractible manifold of dimension $m$. Then the space
  $W\times W\setminus\Delta/\Z_2$ of unordered pairs of points in $W$ is
  homotopy equivalent to $\R P^{m-1}$.
\end{lemma}

\begin{proof}
  We may assume that $n>2$ since otherwise $W$ is homeomorphic to $\R^n$. Let
  $U\subset W$ be a (small) open set homeomorphic to $\R^n$. Consider the
  diagram 
$$\begin{matrix}
 
U\times U\setminus\Delta&\hookrightarrow&U\times U\setminus\Delta\\
\downarrow&&\downarrow\\
(U\times U\setminus\Delta)/\Z_2&\hookrightarrow&(U\times
  U\setminus\Delta)/\Z_2
\end{matrix}
$$ 
Note that $U\times U\setminus\Delta$ fibers over $U$ with fiber
  $U\setminus pt\simeq S^{n-1}$; thus $U\times U\setminus\Delta\simeq S^{n-1}$
  and similarly $W\times W\setminus\Delta\simeq S^{n-1}$; moreover, inclusion
$$
U\times U\setminus\Delta\hookrightarrow W\times W\setminus\Delta
$$
  is a
  homotopy equivalence. Since for $n>2$ the two spaces in the first row of the
  above diagram are simply-connected, it follows that $$(U\times
  U\setminus\Delta)/\Z_2\hookrightarrow (W\times W\setminus\Delta)/\Z_2$$
  induces an isomorphism in homotopy groups, and is therefore a homotopy
  equivalence. 
\end{proof}

In the simplest instance, the lemma below states that the utilities
graph embedded in $\R^3$ links every push-off of itself.

\begin{lemma}[The Linking Lemma] \label{linking} 
  Suppose $W$ is a contractible $(m+1)$-manifold, $K$ is an $m$-obstructor
  complex and $G:K\times [0,\infty) \to W$ is a (continuous) proper map.
  Then there exist two disjoint simplices $\sigma,\tau$ in $K$ such that
  $G(\sigma\times \{0\})\cap G(\tau\times [0,\infty))\neq\emptyset$.
\end{lemma}

\begin{proof} Again we first consider the case $W=\R^{m+1}$.
  Assuming the contrary, consider the homotopy $H_t:\tilde \Sigma\to S^{m}$
  defined by declaring that $H_t(x,y)$ is the class of parallel rays
  containing the ray from $G(x,0)$ through $G(y,t)$. Then $H_0:\tilde\Sigma\to
  S^{m}$ covers a classifying map and therefore has degree 1. Let $B$ be a
  Euclidean ball centered at the origin containing $G(K\times \{0\})$ and
  assume that $t$ is chosen so that $G(K\times \{t\})\cap B=\emptyset$. There
  is a homotopy $L_s$ of $H_t$ defined by setting $L_s(x,y)$ be the
  equivalence class of rays containing the ray from $(1-s)G(x,0)$ through
  $G(y,t)$. Now $L_1$ visibly factors through the projection
  $\pi:\tilde\Sigma\to K$ and therefore by item 3 in the definition of
  obstructor complexes the degree of $L_1$, and
  hence of $H_t$, is 0.  Contradiction.
  
  For the case of a general $W$, replace the definition of $H_t$ by
  $H_t(x,y)=(G(x,0),G(y,t))\in W\times W\setminus \Delta$, and replace the
  ball $B$ by a compact set in which $G(K\times\{0\})$ can be homotoped to a
  point. 
\end{proof}

\subsection{Examples}
The $n$-complexes of van Kampen and of Flores are $(2n)$-obstructor complexes
in our terminology\footnote{Although we will not need the Flores' complexes in what follows.}. 
The collection $\Sigma$ consists of all pairs of disjoint
$n$-simplices. The case of the iterated join of three points can be verified
inductively noting that the three-point set is a 0-obstructor complex and
using the Join Lemma. 
Note that item 3 is vacuous in both examples (for $n>0$).

\begin{lemma}[The Cone Lemma]\label{cone}
  If $K$ is an $m$-obstructor complex, then the cone $cK$ is an
  $(m+1)$-obstructor complex.
\end{lemma}

\begin{proof}
  Let $\Sigma=\Sigma(K)=\{\{\sigma_i,\tau_i\}\}$ be the cycle for $K$. We
  define $\Sigma(cK)$ to have twice as many elements: for every
  $\{\sigma_i,\tau_i\}\in \Sigma$ put $\{c\sigma_i,\tau_i\}$ and
  $\{\sigma_i,c\tau_i\}$ into $\Sigma(cK)$. It is straightforward to check
  items 1 and 3.

  To verify item 2, choose a general position map $f:K\to\R^m$, and let
  $(f,g):K\to \R^m\times\R$ be a perturbation to a general position map. Put
  the cone point high above the hyperplane $\R^m\times\{0\}$ and let $\tilde
  G:cK\to \R^{m+1}$ be the natural extension of $(f,g)$. Then 
$$
|\tilde
  G(c\sigma_i)\cap \tilde G(\tau_i)|+|\tilde G(\sigma_i)\cap \tilde
  G(c\tau_i)|=|f(\sigma_i)\cap f(\tau_i)|  
   $$
  and the claim follows.
\end{proof}

\begin{lemma}[The Join Lemma]\label{join}
  If $K_j$ is an $m_j$-obstructor complex for $j=1,2$ then the join $K_1*K_2$
  is an $(m_1+m_2+2)$-obstructor complex.
\end{lemma}

\begin{proof}
  Let $\Sigma^j=\Sigma(K_j)=\{\{\sigma_i^j,\tau_i^j\}\}$ be the cycle for
  $K_j$, $j=1,2$. We define $\Sigma(K_1*K_2)$ to have $2|\Sigma^1||\Sigma^2|$
  elements: for each $\{\sigma^1_i,\tau^1_i\}\in \Sigma^1$ and
  $\{\sigma^2_l,\tau^2_l\}\in \Sigma^2$ we put the following two pairs in
  $\Sigma(K_1*K_2)$: $\{\sigma^1_i*\sigma^2_l,\tau^1_i*\tau^2_l\}$ and
  $\{\sigma^1_i*\tau^2_l,\tau^1_i*\sigma^2_l\}$. Item 3 is vacuous for
  $\Sigma(K_1*K_2)$ as there are no $(m_1+m_2+2)$-simplices in
  $\Sigma(K_1*K_2)$. 
  
  To verify that $\Sigma(K_1*K_2)$ is a cycle, suppose first that
  $\sigma*\tau$ and $\sigma'*\tau'$ are disjoint simplices of $K_1*K_2$ (with
  $\sigma$ and $\sigma'$ simplices of $K_1$ and $\tau, \tau'$ simplices of
  $K_2$) and that the sum of their dimensions is $m_1+m_2+1$. If
  $\dim(\sigma)+\dim(\sigma')>m_1$ or if $\dim(\tau)+\dim(\tau')>m_2$ then the
  corresponding $(m_1+m_2+1)$-cell is not a face of any $(m_1+m_2+2)$-cells in
  $\Sigma(K_1*K_2)$. So without loss of generality we may assume that
  $\dim(\sigma)+\dim(\sigma')=m_1$ and $\dim(\tau)+\dim(\tau')=m_2-1$. Since
  $\{\tau,\tau'\}$ represents an $(m_2-1)$-cell, item 1 for $K_2$ implies that
  there is an even number $m_2$-cells $\{\tilde\tau_p,\tau'\}$ and
  $\{\tau,\tilde\tau'_q\}$ in $\Sigma^2$ that contain $\{\tau,\tau'\}$. Since
  $\{\sigma*\tilde\tau_p,\sigma'*\tau'\}$ and
  $\{\sigma*\tau,\sigma'*\tilde\tau'_q\}$ are precisely the
  $(m_1+m_2+2)$-cells in $\Sigma(K_1*K_2)$ that contain
  $\{\sigma*\tau,\sigma'*\tau'\}$ the verification of item 1 in this case is
  finished. 
  
  Now suppose that $\sigma*\tau$ and $\sigma'*\tau'$ are disjoint
  simplices of $K_1*K_2$ and that the sum of their dimensions is
  $m_1+m_2+1$. The number of ways of enlarging this cell to an
  $(m_1+m_2+2)$-cell in $\Sigma(K_1*K_2)$ is either 0 (if
  $\{\sigma,\sigma'\}\notin\Sigma^1$) or it equals the number of vertices
  $v\in K_2$ such that $\{\tau,v\}\in\Sigma^2$, which is even by item 3 for
  $K_2$. Thus item 1 is verified for $\Sigma(K_1*K_2)$.
  
  It remains to verify item 2. Let $f_j:K_j\to \R^{m_j}$ be general position
  maps and $I_j$ the total number of intersection points of $f_j$-images of
  unordered pairs of simplices in $\Sigma^j$. View $\R^{m_i}$ as
  $\R^{m_i}\times \{0\}\subset \R^{m_i+1}$ and $\R^{m_1+m_2+2}$ as
  $\R^{m_1+1}\times \R^{m_2+1}$. Perturb $f_j$ to a general position map
  $\tilde G_j=(f_j,g_j):K_j\to \R^{m_j}\times\R=\R^{m_j+1}$ and let
  $G:K_1*K_2\to \R^{m_1+m_2+2}$ be the linear join of $\tilde G_1$ and $\tilde
  G_2$. The number of intersection points of $G$-images of unordered pairs of
  simplices in $\Sigma(K_1*K_2)$ is $I_1I_2$. The details are left to the
  reader.
\end{proof}

\section{The main theorem}

Recall that a (continuous) map $h:A\to B$ is {\it proper} if the preimages of compact sets
are compact. We say that maps $h_1:A_1\to B$ and $h_2:A_2\to B$ into a metric
space $B$ {\it diverge} (from each other) if for every $D>0$ there are compact
sets $C_i\subset A_i$ such that $h_1(A_1\setminus C_1)$ and $h_2(A_2\setminus
C_2)$ are $>D$ apart. If $K$ is a finite complex, we define the open cone
$cone(K)=K\times [0,\infty)/K\times\{0\}$. If $K$ is also an obstructor
complex, we say that a proper map $h:cone(K)\to B$ is {\it expanding} if for
disjoint simplices $\sigma,\tau$ in $K$ the maps $h|cone(\sigma)$ and
$h|cone(\tau)$ diverge. It will also be convenient to make the analogous
definition on the level of 0-skeleta. Triangulate $cone(K)$ so that
$cone(\sigma)$ is a subcomplex whenever $\sigma$ is a simplex of $K$. We say
that a proper map $h:cone(K)^{(0)}\to B$ is {\it expanding} if for all pairs
$\sigma,\tau$ of disjoint simplices in $K$ the restrictions
$h|cone(\sigma)^{(0)}$ and $h|cone(\tau)^{(0)}$ diverge. We also equip
$cone(K)^{(0)}$ with the edge-path metric, so that
a map $h:cone(K)^{(0)}\to B$ is Lipschitz if there is a uniform upper 
bound on the distance between the images of adjacent vertices in $cone(K)^{(0)}$.

Note that if $h:cone(K)\to B$ is a proper expanding map, then there is
$t_0\geq 0$ such that the map $G:K\times [0,\infty)\to B$ defined by
$G(x,t)=h([x,t+t_0])$ satisfies the hypotheses of the Linking Lemma
\ref{linking}, namely $G(\sigma\times\{0\})\cap G(\tau\times
[0,\infty))=\emptyset$ for any two disjoint simplices $\sigma,\tau$ of $K$.

A proper map $h:A\to B$ between proper metric spaces is {\it
  uniformly proper} if there is a proper function $\phi:[0,\infty)\to
[0,\infty)$ such that 
$$
d_B(h(x),h(y))\geq \phi(d_A(x,y))$$
 for all $x,y\in A$. This notion is weaker than the notion of a 
 quasi-isometric embedding,  which would require $\phi$ to be a linear function.

Let $\Gamma$ be a finitely generated group equipped with the
word-metric with respect to some finite generating set. We make the following
definitions.

\begin{definition}\label{obdim}
  The {\it obstructor dimension} $\obdim(\Gamma)$ is defined to be 0 for finite
  groups, 1 for 2-ended groups, and otherwise $m+2$ where $m$ is the largest
  integer such that for some $m$-obstructor complex $K$ and some triangulation
  of the open cone $cone(K)$ as above there exists a proper, Lipschitz,
  expanding map $f:cone(K)^{(0)}\to \Gamma$. If no maximal $m$
  exists we set $\obdim(\Gamma)=\infty$.
\end{definition}

\begin{remark}
  Clearly, one can replace $\Gamma$ in the above definition by any
  quasi-isometric proper metric space. In particular, if $\Gamma$ acts
  cocompactly, properly discontinuously, and isometrically on a proper
  geodesic metric
  space $X$, we can substitute $X$ for $\Gamma$. Moreover, if $\Gamma$ is of
  type $F_{m+1}$ (see e.g., \cite{Brown}) so that $X$ can be chosen to be $m$-connected, then $f$ can
  be extended to a proper, expanding map $\tilde f:cone(K)\to X$ with a
  uniform bound on the diameter of the image of any simplex. One advantage of
  having the (continuous) map defined on the whole cone is that the
  requirement that the map be Lipschitz can be dropped: one can always
  triangulate $cone(K)$ to make the same map Lipschitz.

Note that if $\Gamma$ is infinite and not 2-ended, then we can take $K$ to
consist of 3 points, so $\obdim(\Gamma)\geq 2$.
\end{remark}

\begin{definition}
  The {\it uniformly proper dimension} $\qedim(\Gamma)$ is the
  smallest integer $n$ such that there is a contractible $n$-manifold $W$
  equipped with a proper metric $d_W$ so that there is a Lipschitz,
  uniformly proper map $g:\Gamma\to W$ and so that in addition there is a {\it
    contractibility function} $\rho:(0,\infty)\to (0,\infty)$ such that any
  ball of radius $r$ centered at a point of the image of $g$ is contractible
  in the ball of radius $\rho(r)$ centered at the same point. If no such $n$
  exists we set $\qedim(\Gamma)=\infty$.
\end{definition}

\begin{remark}
  One usually requires of the contractibility function that the statement
  about balls be true regardless of where the center is. If we omit the
  requirement altogether, the invariant would be trivial: every finitely
  generated group admits a uniformly proper map into $[0,\infty)$. Just choose
  an injective map $g:\Gamma\to {\mathbb N}\subset [0,\infty)$. The largest
  metric on $[0,\infty)$ that makes $g$ 1-Lipschitz and makes all $[n,n+1]$
  isometric to a standard closed interval (of length dependent on $n$) is
  proper. Of course, this metric is not a path-metric, but insisting on
  path-metrics would only raise the dimension by 1: For every $\Gamma$ there
  is a proper path-metric on ${\mathbb R}^2$ and a uniformly proper map
  $\Gamma\to\R^2$.
\end{remark}

\begin{definition} 
  The {\it action dimension} $\actdim(\Gamma)$ is the smallest
  integer $n$ such that $\Gamma$ admits a properly discontinuous action on a
  contractible $n$-manifold. If no such $n$ exists, then
  $\actdim(\Gamma)=\infty$.
\end{definition}

Denote also by $\gdim(\Gamma)$ the {\it geometric dimension} of $\Gamma$, i.e., 
the minimal $n$ such that $\Gamma$ admits a properly discontinuous action on a
contractible $n$-complex. Recall that for virtually torsion-free groups
$\Gamma$, $\gdim(\Gamma)$ is conjectured to be equal to the virtual
cohomological dimension $\vcdim$ of $\Gamma$ and that the only potential
counterexamples would have $\gdim=3$ and $\vcdim=2$ (see \cite{Brown}).

We note that $\qedim(\Gamma)\leq \actdim(\Gamma)$ by choosing a proper
invariant metric on $W$ and taking an orbit of the action, and that for
torsion-free groups $\Gamma$ we have $\actdim(\Gamma)\leq 2\cdot
\gdim(\Gamma)$ by the Stallings theorem cited in the
introduction. Alternatively, we could find a $(2n)$-dimensional thickenning of
an $n$-complex by immersing it in $\R^{2n}$ and taking a regular
neighborhood. The inequality $\actdim(\Gamma)\leq 2\cdot
\gdim(\Gamma)$ is false for groups with torsion; indeed, the free product
$A_5*A_5$ acts properly discontinuously on a tree but not on the plane.
On the other hand, $\Gamma=A_5* A_5$ contains a free subgroup $\Gamma'$ of finite 
index, hence $2=\actdim(\Gamma')< \actdim(\Gamma)$. 

The main
theorem in this note is:

\begin{thm} \label{main2} 
  $\obdim(\Gamma)\leq \qedim(\Gamma)$.
\end{thm}

\begin{proof} The special cases when $\obdim(\Gamma)\leq 1$ are clear. 
  Let $K$ be an $m$-obstructor complex and $f:cone(K)^{(0)}\to \Gamma$ a
  proper, Lipschitz, expanding map. Let $W$ be a contractible manifold with a
  proper metric and $g:\Gamma\to W$ a uniformly proper Lipschitz map
  satisfying the contractibility function requirement.  Consider the
  composition $gf:cone(K)^{(0)}\to W$. Now extend $gf$ inductively over the
  skeleta of $cone(K)$ to get a map $G:cone(K)\to W$. Using the
  contractibility function, we can arrange that the diameter of the image of
  each simplex of $cone(K)$ is uniformly bounded. It follows that $G$ is a
  proper expanding map, and therefore $n\geq m+2$ by the Linking Lemma
  \ref{linking}.
\end{proof}

This theorem immediatately implies the following chain of
inequalities (with the last inequality only for torsion-free groups):

\begin{equation}\label{in}
\obdim(\Gamma)\leq \qedim(\Gamma)\leq \actdim(\Gamma)\leq 2
\gdim(\Gamma)
\end{equation}

The second inequality can be strict. The Baumslag-Solitar group
$$
B=\<x,t|xt=t^2x\>$$
is not a 3-manifold group and so $\actdim(B)=2\cdot
\gdim(B)=4$. On the other hand, $\obdim(B)=\qedim(B)=3$. The group $B$ admits a
uniformly proper map into ${\mathbb H}^3$ and the universal cover of the
presentation 2-complex admits an expanding proper homotopy of the tripod,
which is a 1-obstructor complex. All three invariants in (\ref{in}) 
are monotone, in the
sense that if $\Gamma'$ is a finitely generated subgroup of $\Gamma$, then
$anydim(\Gamma')\leq anydim(\Gamma)$.  We also note that 
both $\obdim$ and $\qedim$ are 
invariant under quasi-isometries. This is not the case for $\actdim$ 
(even for torsion-free groups) as there
are examples of torsion-free groups that are not 3-manifold groups but contain
3-manifold groups as finite index subgroups \cite{kk:duality}.

\begin{lemma}
$$\obdim(\Gamma_1\times\Gamma_2)\geq \obdim(\Gamma_1)+ \obdim(\Gamma_2)$$
while
$$\qedim(\Gamma_1\times\Gamma_2)\leq \qedim(\Gamma_1)+ \qedim(\Gamma_2)$$
and
$$
\actdim(\Gamma_1\times\Gamma_2)\leq \actdim(\Gamma_1)+ \actdim(\Gamma_2).$$
\end{lemma}
\begin{proof} The latter two statements are obvious, while the first one follows from 
the Join Lemma \ref{join}. The product $cone(K_1)\times cone(K_2)$ can naturally
be viewed as $cone(K_1*K_2)$ and the product map into $\Gamma_1\times\Gamma_2$
satisfies the requirements. (If one of the two groups is 2-ended, use the Cone
Lemma instead.) 
\end{proof} 

\begin{cor}
In particular, we see that for $\Gamma=F_2^n$ all three
invariants $\obdim, updim$ and $\actdim$ are $2n$ and the 
inequalities in the chain (1) are equalities.
\end{cor}

If $\Gamma$ has a reasonable boundary, it may be easier to compute
$\obdim(\Gamma)$. The following definition is taken from \cite{mb:boundary}.

\begin{definition} \label{zstructure}
Let $\Gamma$ be a group. A $\cal Z$-{\it structure} on $\Gamma$
is a pair $(\tilde X,Z)$ of spaces satisfying the following four
axioms.
\begin{itemize}
\item $\tilde X$ is a Euclidean retract.
\item $Z$ is a $Z$-set in $\tilde X$.
\item $X=\tilde X\setminus Z$ admits a covering space
action of $\Gamma$ with compact quotient.
\item The collection of translates of a compact set in $X$ forms a
null-sequence in $\tilde X$, i.e., for every open cover $\cal U$ of $\tilde X$
all but finitely many translates are $\cal U$-small.
\end{itemize}
A space $Z$ is a {\it boundary} of $\Gamma$ if there is a $\cal
Z$-structure $(\tilde X,Z)$ on $\Gamma$.   
\end{definition}

For example, torsion-free hyperbolic groups and $CAT(0)$ groups\footnote{I.e. 
groups which admit discrete cocompact isometric action on a $CAT(0)$-space.} 
admit a boundary. However, unlike for hyperbolic groups, boundary of a $CAT(0)$-group 
$G$ is not uniquely determined by $G$ (up to a homeomorphism) \cite{kleiner-croke}. 

\begin{cor} \label{boundary}
  Suppose $Z$ is a boundary of $\Gamma$ and $f:K\to Z$ is a map from an
  $m$-obstructor complex that sends disjoint simplices disjointly (e.g. $f$
  could be an embedding). Then $\obdim(\Gamma)\geq m+2$.\qed
\end{cor}

\begin{proof}
  Let $X,\tilde X$ be as in the definition. Since $Z$ is a $Z$-set in $\tilde
  X$, there is a homotopy $H:K\times [0,1]\to \tilde X$ with $H(x,0)=*$,
  $H(x,1)=f(x)$ and $H(K\times (0,1])\cap Z=\emptyset$. Restricting to
  $K\times [0,1)$ and reparametrizing yields an expanding map $cone(K)\to X$.
\end{proof}

It is convenient to introduce the notation 
$$
``K\subset\partial\Gamma"
$$ 
to mean that there is a proper expanding Lipschitz map $cone(K)^{(0)}\to\Gamma$
as in the definition of $\obdim$. The above corollary implies
$$
K\subset\partial \Gamma\Rightarrow ``K\subset\partial \Gamma"
$$

\begin{example}
  The $n$-fold join of Cantor sets is a boundary of $F_2^n$ and it contains
  van Kampen's $(2n-2)$-obstructor complex. Thus $\obdim(F_2^n)=2n$ and all
  inequalities in the chain (\ref{in})  are equalities.
\end{example}

\begin{remark}
It seems to be believed by the experts that there are $n$-dimen\-sio\-nal
torsion-free hyperbolic groups $\Gamma$ with boundary the Menger universal
$(n-1)$-dimensional compactum. For such a group all inequalities would be
equalities as well, but no such examples of hyperbolic groups are known except
for small $n$.
\end{remark}

Somewhat more generally, consider a group $G$ acting discretely 
isometrically on a CAT(0)-space $X$ with the ideal boundary 
$D=\partial_{\infty} X$. (We do not assume that this action is 
cocompact.) 
Pick a base point $x\in X$ and $C\in \R_+$. 
The $C$-cone limit set $\Lambda_C(G)$ of $G$  consists of points 
$\xi\in D$ such that for the geodesic ray $\rho$ in $X$ emanating 
from $x$ and representing $\xi$, there exists 
an infinite sequence $g_n\in G$ such that $d(g_n x, \rho) \le C$. 
The arguments from Corollary \ref{boundary} imply

\begin{cor}\label{CC}
If for some $C$, $\Lambda_C(G)$ contains an $m$-obstructor complex, 
then $\obdim(G)\ge m$. 
\end{cor}

\section{Short exact sequences}

We now investigate the obstructor dimension of a group $G$ that fits in a
short exact sequence
$$
1\to H\to G\overset{\pi}\to Q\to 1
$$
where all groups are finitely generated. The natural guess is that
\begin{equation}\label{seq}
\obdim G\geq \obdim H+\obdim Q
\end{equation}
and this is what we prove under certain technical assumptions on $\pi$ (admits
a Lipschitz section) and $H$ (weakly convex). All groups are equipped with
word metrics.  We note that {\em some} restrictions are clearly neccessary for
(\ref{seq}) to hold. For instance, Rips in \cite{Rips} constructs examples of
2-dimensional hyperbolic groups $G$ which admit epimomorphisms $G\to Q$ (where
$Q$ is the prescribed finitely presented group) so that the kernel $H$ is
finitely generated (and is neither finite nor 2-ended).  Note that $G$ can be
assumed to have Menger curve boundary \cite{kk}.  Then $obdim(G)=4$, $\obdim
H\ge 2$ and $Q$ can be chosen to have $\obdim(Q)$ as large as one likes.
See also Example \ref{CE}. 

\begin{definition}
  We say that a finitely generated group $\G$ is
  {\it weakly convex} if there is a collection of (discontinuous, of course)
  paths $\{\phi_{z,w}:[0,1]\to \G\}_{z,w\in\G}$ and a constant $M>0$
  satisfying the following properties:
\begin{enumerate}
\item $\phi_{z,w}(0)=z$ and $\phi_{z,w}(1)=w$.
\item There is a function $\gamma:[0,\infty)\to [0,\infty)$ such that 
$$d(z,w)\leq R\Longrightarrow \diam(Im(\phi_{z,w}))\leq \gamma(R).$$
\item For all $z,w\in\G$ there is $\epsilon>0$ such that $\phi_{z,w}$ sends
  subintervals of length $<\epsilon$ to sets of diameter $<M$.
\item If $d(z,z')\leq 1$ and $d(w,w')\leq 1$ then for all $t\in [0,1]$ 
$$d(\phi_{z,w}(t),\phi_{z',w'}(t))\leq M.$$
\end{enumerate}
\end{definition}

\begin{remark}
  The paths are to be thought of as being piecewise constant. We could avoid
  talking about discontinuous functions by requiring that they be defined only
  on the rationals in $[0,1]$. It is more standard to think of paths in $\G$
  as eventually constant 1-Lipschitz functions defined on non-negative
  integers; however, for what follows it is important that all paths be
  defined on the same bounded set. It is 
  possible to reparametrize such paths by ``constant speed'' paths defined on
  $[0,1]$.  The collection of paths as above is usually called a ``combing''
  (except for the domain being $[0,1]$).  Condition 2 is then a weak version
  of the requirement that the combing be quasi-geodesic and it follows
  automatically if the combing is equivariant (i.e., $\phi_{gz,gw}=L_g\circ
  \phi_{z,w}$, where $L_g:\Gamma\to\G$ denotes left translation by $g$).
  Condition 3 is the replacement of the 1-Lipschitz requirement. Condition 4
  is the ``Fellow Traveller'' property.
  
  If $\Gamma$ and $\G'$ are quasi-isometric and one is weakly convex, so is
  the other. Hyperbolic groups, $CAT(0)$ groups, and semi-hyperbolic groups
  \cite{ab:semihyperbolic} are weakly convex.

\end{remark}

We can regard the given paths in the definition of weak convexity as a recipe
for extending maps into $\G$ defined on the (ordered) vertices of a 1-simplex
to the whole 1-simplex. It is easy to see that one can similarly extend maps
defined on the vertices of an $n$-simplex for any $n>0$, with the constant
$M=M(n)$ above depending on $n$. By
$$\Delta^n=\{(t_0,t_1,\cdots,t_n)\in\R^{n+1}|t_i\geq
0,t_0+t_1+\cdots+t_n=1\}$$ we denote the standard $n$-simplex, and by
$I_{n,k}$ the standard face inclusion $\Delta^{n-1}\hookrightarrow \Delta^n$
onto the face $t_k=0$ given by 
$$I_{n,k}(t_0,t_1,\cdots,t_{n-1})=(t_0,t_1,\cdots,0,\cdots,t_{n-1}).$$

\begin{prop}\label{extend}
Let $\G$ be a weakly convex group. Then for every $n>0$ there is a constant
$M(n)$ and for every $(n+1)$-tuple $(z_0,z_1,\cdots,z_n)\in\G^{n+1}$ there
is a function $\phi_{z_0,z_1,\cdots,z_n}:\Delta^n\to\G$ such that
\begin{itemize}
\item $\phi_{z_0,z_1,\cdots,z_n}(v_k)=z_k$ where $v_k\in\Delta^n$ is the
  vertex with $t_k=1$.
\item There is a function $\gamma_n:[0,\infty)\to [0,\infty)$ such that
$$
d(z_i,z_j)\leq T\text{ for all }i,j\Longrightarrow
\diam(Im(\phi_{z_0,z_1,\cdots,z_n}))\leq \gamma_n(T)
$$
\item For all $z_0,z_1,\cdots,z_n$ there is $\epsilon>0$ such that the
  $\phi_{z_0,z_1,\cdots,z_n}$-images of sets of diameter $<\epsilon$ have
  diameter $<M(n)$.
\item If $d(z_i,w_i)\leq 1$ then
  $$d(\phi_{z_0,z_1,\cdots,z_n}(t),\phi_{w_0,w_1,\cdots,w_n}(t))\leq M(n)$$
\item $\phi_{z_0,z_1,\cdots,z_n}\circ I_{n,k}=\phi_{z_0,z_1,\cdots,\hat
    z_k,\cdots,z_n}$
\end{itemize}
\end{prop}

\begin{proof}[Proof (sketch)]
Functions $\phi_{z_0,z_1,\cdots,z_n}$ are constructed by induction on $n$,
with the case $n=1$ being the definition. The inductive step consists of
defining $\phi_{z_0,z_1,\cdots,z_n}$ on the boundary of $\Delta^n$ so that the
last item above holds and then extending to the interior by coning off from
the first vertex. More precisely, if $\psi:[0,1]\to\Delta^n$ is a linear map
with $\psi(0)=v_0$ and $\psi(1)$ belongs to the face with $t_0=0$ then 
$$\phi_{z_0,z_1,\cdots,z_n}(\psi(t))=\phi_{z_0,w}(t)$$ where $w\in\G$ is the
image of $\psi(1)$ under the (partially defined)
$\phi_{z_0,z_1,\cdots,z_n}$. Checking the properties listed above is
straightforward (by construction, the restriction of the function to a face
containing $v_0$ is already the cone on the opposite face).
\end{proof}

\begin{thm}\label{semi}
  Let $$1\to H\to G\overset{\pi}\to Q\to 1$$
  be a short exact sequence of
  finitely generated groups. Suppose that $H$ is weakly convex and that $\pi$
  admits a Lipschitz section $s:Q\to G$. Then $$\obdim G\geq \obdim H+\obdim
  Q$$ 
\end{thm}

\begin{proof}
If $H$ (or $Q$) is finite, then $G$ is quasi-isometric to $Q$ (or $H$) and
  equality holds. If $H$ (or $Q$) is 2-ended, we can use $K_H=point$ (or
  $K_Q=point$) in the proof below and appeal to the Cone Lemma. If both $H$
  and $Q$ are 2-ended, then $G$ is virtually $\Z\times\Z$ and thus
  $\obdim(G)=2$, $\obdim(H)=\obdim(Q)=1$, so equality again holds.

  Let $\alpha:cone(K_H)^{(0)}\to H$ and $\beta:cone(K_Q)^{(0)}\to Q$ be proper
  Lipschitz expanding maps defined on the vertices of a fine triangulation of
  the cones on obstructor complexes $K_H$ and $K_Q$. Define
  $$f:cone(K_H*K_Q)=cone(K_H)^{(0)}\times cone(K_Q)^{(0)}\to G$$
  by
$$f(x,y)=\alpha(x)\cdot s\beta(y).$$ 

{\bf Claim 1.} $f$ is a proper map.

Indeed, let $(x_i,y_i)$ be a sequence in $cone(K_H)^{(0)}\times
cone(K_Q)^{(0)}$ leaving every finite set. If the sequence $\pi
f(x_i,y_i)=\beta(y_i)\in Q$ leaves every finite set, the same is true for
$f(x_i,y_i)\in G$. Otherwise, after passing to a subsequence, we may assume
that the sequence $\pi f(x_i,y_i)=\beta(y_i)\in Q$ stays in a finite set
$D\subset Q$. Then $s\beta(y_i)$ stays in the finite set $s(D)$. Since $\beta$
is a proper map, the sequence $y_i\in cone(K_Q)^{(0)}$ stays in a finite set,
and thus the sequence $x_i\in cone(K_H)^{(0)}$ leaves every finite set. Since
$\alpha$ is a proper map, we see that the sequence
$f(x_i,y_i)=\alpha(x_i)\cdot s\beta(y_i)$ leaves every finite set.

{\bf Claim 2.} If $\sigma=\sigma_H*\sigma_Q$ and $\tau=\tau_H*\tau_Q$ are
disjoint simplices of $K_H*K_Q$, then $f|cone(\sigma)^{(0)}$ and
$f|cone(\tau)^{(0)}$ 
diverge. 

Indeed, let $(x_i,y_i)$ and $(x_i',y_i')$ be sequences in
$cone(\sigma_H)^{(0)}\times cone(\sigma_Q)^{(0)}$ and
$cone(\tau_H)^{(0)}\times cone(\tau_Q)^{(0)}$ respectively, leaving every
finite set. Note that $\pi$ is a Lipschitz map, so if one of two sequences
$\pi f(x_i,y_i)=\beta(y_i)$ and $\pi f(x_i',y_i')=\beta(y_i')$ leaves every
finite set in $Q$, then $d_Q(\beta(y_i),\beta(y_i'))\to\infty$ (since
$\beta|cone(\sigma_Q)^{(0)}$ and $\beta|cone(\tau_Q)^{(0)}$ diverge) and
consequently $$d_G(f(x_i,y_i),f(x_i',y_i'))\to\infty.$$ Now assume that both
sequences $\beta(y_i)$ and $\beta(y_i')$ are contained in a fixed finite set
$D\subset Q$. Then we have
\begin{equation*}
\begin{split}
d_G(f(x_i,y_i),f(x_i',y_i'))=d_G(\alpha(x_i)\cdot
s\beta(y_i),\alpha(x_i')\cdot s\beta(y_i'))=\\
d_G(1,s\beta(y_i)^{-1}\alpha(x_i)^{-1}\alpha(x_i')s\beta(y_i'))
\end{split}\end{equation*}
Since $s\beta(y_i)$ and $s\beta(y_i')$ stay in a finite set and
$$d_H(1,\alpha(x_i)^{-1}\alpha(x_i'))=d_H(\alpha(x_i),\alpha(x_i'))\to\infty$$
it follows that $d_G(1,\alpha(x_i)^{-1}\alpha(x_i'))\to\infty$ and
$$d_G(1,s\beta(y_i)^{-1}\alpha(x_i)^{-1}\alpha(x_i')s\beta(y_i'))\to\infty,$$
and the claim is proved.

The remaining problem is that $f$ is not Lipschitz. 

{\bf Claim 3.} The restriction of $f$ to $\{q\}\times cone(K_Q)^{(0)}$ is
Lipschitz with the Lipschitz constant independent of $q$.

Indeed, let $x,y$ be two adjacent vertices in $cone(K_Q)^{(0)}$. 
$$d_G(f(q,x),f(q,y))=d_G(\alpha(q)\cdot s\beta(x),\alpha(q)\cdot s\beta(y))=
d_G(s\beta(x),s\beta(y))$$
and the claim follows from the assumption that $s$
is Lipschitz.

For every $p\in cone(K_Q)^{(0)}$ let $f_p$ denote the restriction of $f$ to
the slice $cone(K_H)^{(0)}\times \{p\}$. Recall that $L_{s\beta(p)}$ is the
left translation by $s\beta(p)$ and it induces an isometry between
$H=\pi^{-1}(1)$ (with the $G$-metric) and $\pi^{-1}(\beta(p))$.

{\bf Claim 4.} 
$L_{s\beta(p)}^{-1}f_p:cone(K_H)^{(0)}\times \{p\}\to H$
is Lipschitz with respect to the word-metric on $H$ (but the Lipschitz
constant depends on $p$). In particular, $f_p:cone(K_H)^{(0)}\times \{p\}\to
G$ is Lipschitz.

Indeed, $L_{s\beta(p)}^{-1}f_p(x,p)= s\beta(p)^{-1}\cdot\alpha(x)\cdot
s\beta(p)$ which is Lipschitz.

We next order all vertices of $cone(K_H)$ and then extend (simplex-by-simplex)
$L_{s\beta(p)}^{-1}f_p$ for each $p$ to the map $\tilde
F_p:cone(K_H)\times\{p\}\to H$ using the weak convexity of $H$ and Proposition
\ref{extend}. Then define
$$
\tilde f_p=L_{s(p)}\tilde F_p:cone(K_H)\times\{p\}\to \pi^{-1}(p).$$
Let
$$
\tilde f:cone(K_H)\times cone(K_Q)^{(0)}\to G$$
be defined as $\tilde f_p$ on each $cone(K_H)\times\{p\}$.

We now note that for $n=\dim cone(K_H)$ and for $M=M(n)$ from Proposition
\ref{extend} we have that:
\begin{itemize}
\item For each $p\in cone(K_H)^{(0)}$ there is $\epsilon(p)>0$ so that sets of
  diameter $<\epsilon(p)$ in a simplex of $cone(K_H)\times \{p\}$ are sent by
  $\tilde f$ to sets of diameter $<M$.
\item If $p,p'$ are adjacent vertices in $cone(K_Q)$ then $\tilde f(x,p)$ and
  $\tilde f(x,p')$ are $kM$-close for any $x\in cone(K_H)$, where $k$ is a
  Lipschitz constant for $s\beta$.
\end{itemize}
For each $p\in cone(K_Q)^{(0)}$ choose a positive integer $m(p)$ so that the
simplices of the $m(p)^{th}$ barycentric subdivision of $\sigma\times\{p'\}$
have diameter $<\epsilon(p')$ for all vertices $p'$ at distance $\leq 1$ from
$p$. 

We now define a triangulation of $cone(K_H)\times cone(K_Q)$. Start with a
decomposition into cells of the form $\sigma\times \tau$ where $\tau$ is a
simplex of $cone(K_Q)$ and $\sigma$ is a simplex of the $k(\tau)^{th}$
barycentric subdivision of a simplex of $cone(K_H)$ with $k(\tau)=\min
\{m(p), p\in\tau^{(0)}\}$. Now triangulate each such cell inductively on the
dimension so that the vertex set of the triangulation is precisely
\begin{equation*}
\begin{split}
\{(v,p)|p\in cone(K_Q)^{(0)},v
\text{ is a vertex of the\ }m(p)^{th}\\ \text{barycentric
 subdivision of
  }cone(K_H)\}\end{split}\end{equation*}
The restriction of $\tilde f$ to the vertex set is now a Lipschitz
function. 

{\bf Claim 5.} This restriction is still proper and expanding.

The proof closely follows proofs of Claims 1 and 2. If the sequence
$\beta(y_i)$ (resp. one of the sequences $\beta(y_i)$ or
$\beta(y_i')$) leaves every finite set, the proof is exactly the same
as in Claims 1 and 2. Otherwise, without loss of generality, the
sequence $\tilde f(x_i,y_i)$ (resp. sequences $\tilde f(x_i,y_i)$ and
$\tilde f(x_i',y_i')$) belong to finitely many slices of the form
$cone(K_H)\times\{p\}$ and therefore lie a bounded distance away from
sequences considered in Claims 1 and 2 coming from the vertices of the
original triangulation. The proof follows.
\end{proof}

\begin{cor}\label{semi2}
If $G=H\rtimes Q$ with $H$ and $Q$ finitely generated and $H$ weakly convex,
then $\obdim G\geq \obdim H+\obdim Q$.
\end{cor}

\begin{example}\label{CE}
  Let $G= F^n_2$. Define $\phi:G\to\Z$ by sending the basis elements of each
  factor to $1\in\Z$. Let $H=Ker(\phi)$. It is easy to see that $H$ contains a
  copy of $F^n_2$ and thus $\obdim(G)=\obdim(H)=2n$.  Therefore $G= H\rtimes
  \Z$, but $\obdim(H)+\obdim(\Z)=2n+1 > \obdim(G)$.  It follows that $H$ is
  not weakly convex, and one knows \cite{bieri} that $H$ is of type $F_{n-1}$
  (in particular, it is finitely generated for $n\ge 2$, finitely presented
  for $n\ge 3$, etc.).
\end{example}

\begin{remark}
Example \ref{CE} shows that in the above corollary 
weak convexity of $H$ is a necessary assumption. 
\end{remark}

We now apply the above theorem to the group $Out(F_n)$ of outer automorphisms
of $F_n=\<x_1,\cdots,x_n\>$, the free group of rank $n$.  Recall
\cite{cv:moduli} that the virtual cohomological dimension of $Out(F_n)$ is
$2n-3$ ($n>1$). It follows from the Stallings theorem that
$\obdim(Out(F_n))\leq 4n-6$. We note that equality holds, since $Out(F_n)$
contains as a subgroup a group of the form $F_2^{2n-4}\rtimes F_2$. Indeed,
choose $F_2<Aut(F_2)$ that injects into $Out(F_2)$ and let it act diagonally
on $F_2^{2n-4}$. The corresponding semi-direct product is realized by the
subgroup of $Aut(F_n)$ that injects into $Out(F_n)$ as follows. Send an element
\begin{eqnarray*}
u=(w_3,v_3,...,w_{n},v_{n}, \alpha)\in F_2^{2n-4}\rtimes F_2, 
w_i, v_i\in F_2=\<x_1,x_2\>,\\
 \alpha\in F_2\subset Aut(F_2)
\end{eqnarray*}
to the automorphism $\phi_u$ of $F_n$ which acts on $\<x_1,x_2\>$ as the automorphism  
$\alpha$ and maps $x_i$ ($i\ge 3$) to $w_i x_iv_i^{-1}$. The reader will verify 
that $u\mapsto \phi_u$ is indeed a monomorphism  $F_2^{2n-4}\rtimes F_2$ to 
$Aut(F_n)$ whose image projects injectively to $Out(F_n)$.

  We conjecture that $\obdim$ for the mapping class group of
  a closed oriented surface of genus $g$ is $6g-6$, the dimension of the
  associated Teichm\" uller space. This is analogous to Theorem \ref{2} stated
  in the Introduction.

Let $B_n$ denote the braid group on $n$ strands. As evidence for the
conjecture in the previous paragraph, we note that $\obdim B_n=2n-3$ ($n\geq
2$). Denote by $P_n$ the pure braid group on $n$ strands. Then
$P_n=F_{n-1}\rtimes P_{n-1}$ so the statement that $\obdim P_n=2n-3$ follows
by induction from Corollary \ref{semi2}.

\section{Questions}

We conclude this note with the following questions about the invariants.

\noindent 1. Is $\obdim\Gamma=\qedim\Gamma$ for all $\Gamma$? The answer is
probably negative as stated, but the question should be interpreted liberally:
Is $\qedim\Gamma$ detected homologically? (Compare this with a
theorem of Kuratowsky and Claytor \cite{Claytor} that a 1-dimensional 
continuum without global cut-points is planar provided that 
 it contains neither the complete graph on five vertices, nor the 
``utility" graph.)

\noindent 2. 
Suppose $M_i$ is a compact aspherical $n_i$-manifold with all boundary
components aspherical and incompressible, $i=1,...,k$. 
If $M_i$ is not homotopy equivalent
to an $(n_i-1)$-manifold, and if $G=\pi_1(M_1)\times\cdots\times\pi_1(M_k)$,
is $\actdim G=n_1+\cdots+n_k$?

\noindent 3. Is the assumption of the existence of a Lipschitz section in
Theorem \ref{semi} necessary if in addition $H$, $G$, $Q$ have finite type
(i.e., have finite Eilenberg-MacLane spaces)? 

\noindent 
4. What is  $\actdim(\Gamma)$ for uniform/nonuniform S-arithmetic groups?
   Every such group acts on the product of symmetric spaces and buildings. A
   natural guess is that the answer equals the sum of dimensions of the
   symmetric spaces plus twice the sum of dimensions of the buildings.

\noindent 5. Are there groups of finite type (i.e., groups of type $FP$) 
which are not weakly convex? For instance, are the fundamental groups 
of 3-dimensional Nil-manifolds weakly convex?

\bigskip
{\bf Acknowledgements.} The first two authors were supported by NSF
grants DMS-96-26633 and DMS-99-71404; the third author was supported
by a Sloan Foundation Fellowship, and NSF grants DMS-96-26911 and
DMS-99-72047. The second author is also grateful to the Max--Plank
Instutute (Bonn) for its support.

%\bibliography{$HOME/lib/ref}

\begin{thebibliography}{Whi44b}

\bibitem[AB95]{ab:semihyperbolic}
J.~M. Alonso and M.~R. Bridson, \emph{Semihyperbolic groups}, Proc.
  London Math. Soc. (3) \textbf{70} (1995), no.~1, 56--114.

\bibitem[Bes96]{mb:boundary}
M. Bestvina, \emph{Local homology properties of boundaries of groups},
  Michigan Math. J. \textbf{43} (1996), no.~1, 123--139.

\bibitem[BF]{bf:cusp}
M. Bestvina and M. Feighn, \emph{Proper actions of lattices on
  contractible manifolds}, preprint, 2000.

\bibitem[B76]{bieri}
R. Bieri, \emph{Normal subgroups in duality groups and in groups of
cohomological dimension $2$}, J. Pure Appl. Algebra,
\textbf{1} (1976), 35--51.

\bibitem[Br82]{Brown}
 K.S.~Brown, \emph{Cohomology of Groups},
 Graduate Texts in Math., vol. \textbf{87},
Springer Verlag, 1982. 



\bibitem[Cl24]{Claytor}
S.~Claytor, \emph{Topological immersion of {P}eanian 
continua in the spherical surface}, Annals of Math., 
\textbf{35} (1924), 809--835. 



\bibitem[CK00]{kleiner-croke}
C.~B. Croke and B.~Kleiner, \emph{Spaces with nonpositive curvature
  and their ideal boundaries}, Topology \textbf{39} (2000), no.~3, 549--556.

\bibitem[Coo79]{cooke:thickenings}
G. Cooke, \emph{Thickenings of {C}{W} complexes of the form ${S}\sp{m}\cup
  \sb{\alpha }e\sp{n}$}, Trans. Amer. Math. Soc. \textbf{247} (1979), 177--210.

\bibitem[CV86]{cv:moduli}
M.~Culler and K.~Vogtmann, \emph{Moduli of graphs and automorphisms of free
  groups}, Invent. Math. \textbf{84} (1986), 91--119.

\bibitem[DR93]{dr:stallings}
A.~N. Drani{\v{s}}nikov and D.~Repov{\v{s}}, \emph{Embeddings up to homotopy
  type in {E}uclidean space}, Bull. Austral. Math. Soc. \textbf{47} (1993),
  no.~1, 145--148.

\bibitem[Flo35]{flores:embedding}
A.~Flores, \emph{\"{U}ber $n$-dimensionale {K}omplexe die im ${R}_{2n+1}$
  absolut selbstverschlungen sind}, Ergebnisse eines math. Koll. \textbf{6}
  (1935), 4--6.

\bibitem[FKT]{Freedman-Krushkal-Teichner}
M.~Freedman, V.~Krushkal and P.~Teichner,
\emph{van {K}ampen's embedding obstruction is incomplete for
             $2$-complexes in ${\bf {R}}\sp 4$},
 Math. Res. Letters, \textbf{1} (1994), 167--176. 


\bibitem[KK]{kk:duality}
M.~Kapovich and B.~Kleiner, \emph{Coarse {A}lexander duality and duality
  groups}, preprint, 1999.

\bibitem[KK00]{kk}
M.~Kapovich and B.~Kleiner, \emph{Hyperbolic groups with low-di\-men\-sional
  boundaries}, Ann. Ec. Norm. Sup. Paris (2000), to appear.

\bibitem[Kr00]{Krushkal}
V.~Krushkal, \emph{Embedding obstructions and $4$-dimensional thickenings of
             $2$-complexes}, Proc. Amer. Math. Soc. (2000), to appear. 


\bibitem[Rip82]{Rips}
E.~Rips, \emph{Subgroups of small cancellation groups}, Bull. London Math. Soc.
  \textbf{14} (1982), no.~1, 45--47.

\bibitem[Sha57]{shapiro:embedding}
A. Shapiro, \emph{Obstructions to the imbedding of a complex in a euclidean
  space. {I}. {T}he first obstruction}, Ann. of Math. (2) \textbf{66} (1957),
  256--269.

\bibitem[vK33]{vk:embedding}
E.~R. van Kampen, \emph{Komplexe in euklidischen {R}\" aumen}, Abh. Math. Sem.
  Univ. Hamburg \textbf{9} (1933), 72--78 and 152--153.

\bibitem[Whi44a]{whitney:embeddings}
H. Whitney, \emph{The self-intersections of a smooth $n$-manifold in
  $2n$-space}, Ann. of Math. (2) \textbf{45} (1944), 220--246.

\bibitem[Whi44b]{whitney:immersions}
H. Whitney, \emph{The singularities of a smooth $n$-manifold in
  $(2n-1)$-space}, Ann. of Math. (2) \textbf{45} (1944), 247--293.

\bibitem[Wt65]{wu:embedding}
Wu~Wen-ts{\"u}n, \emph{A theory of imbedding, immersion, and isotopy of
  polytopes in a euclidean space}, Science Press, Peking, 1965.

\end{thebibliography}

\providecommand{\bysame}{\leavevmode\hbox to3em{\hrulefill}\thinspace}

\end{document}